\documentclass[11pt]{article}
\usepackage{graphicx}
\usepackage{amsmath,amsfonts}
\newcommand{\beq}{\begin{equation}}
\newcommand{\eeq}{\end{equation}}
\newcommand{\defeq}{\stackrel{\mbox{\tiny def}}{=}}

\newcommand{\diag}{\mathrm{diag}}
\newcommand{\la}{\lambda}
\newtheorem{exa}{Example}
\newenvironment{example}{\begin{exa}\rm }{\end{exa}}
\newtheorem{rem}{Remark}
\newenvironment{remark}{\begin{rem}\rm }{\end{rem}}
\newcommand{\nn}{\nonumber}

\title{On the spectra of generalized Fibonacci and Fibonacci-like operators}
\author{Ivan Slapni\v{c}ar\thanks{University of Split, Faculty of
    Electrical Engineering, Mechanical Engineering and Naval
    Architecture, R. Boskovica b.b., 21000 Split, Croatia, e-mail:
    ivan.slapnicar@fesb.hr. Author acknowledges the grant number
    023-0372783-1289 of the Ministry of Science, Education
    and Sports of the Republic of Croatia and the grant FP7 People IEF
    ``MATLAN'' of the European Commission.}}

\begin{document}
\maketitle

\begin{abstract}
We analyze the spectra of generalized Fibonacci and Fibonacci-like operators in
Banach space $l^1$. Some of the results have application in population dynamics.
\end{abstract}

\section{Introduction and preliminaries}

Let $l^1$ denote the Banach space of all real sequences 
$x\defeq (x_1,x_2,x_3,\cdots)$ such that $\|x\|_1\defeq\sum |x_k|<
\infty$.
Let $H:l^1\to l^1$ be a linear operator on $l^1$.
The resolvent set of $H$, $\rho(H)$ is the set of all complex numbers $\lambda$
such that the operator $\lambda I-H$ has
a bounded inverse, where $I:l^1\to l^1$ is the identity operator. The set
$\sigma(H)\defeq \mathbb{C}\setminus \rho(H)$ is the
spectrum of  $H$.
The spectrum is further subdivided into three mutually disjoint parts,
the point spectrum $\sigma_p(H)$, the continuous spectrum
$\sigma_c(H)$ and the residual spectrum $\sigma_r(H)$.
The point spectrum is the set of all  $\lambda\in\mathbb{C} $ such that 
$\lambda I-H$ has no inverse. As in the finite dimensional case, such 
$\lambda$ are also called eigenvalues and the corresponding non-zero
vectors $x\in l^1$, such
that $(\lambda I-H)x=0$ are called eigenvectors.
The continuous spectrum is the set of all $\lambda$ not in
$\rho(H)$ or $\sigma_p(H)$ for which the range of $\lambda
I-H$ is dense in $l^1$. The residual spectrum is
the set of all $\lambda$ in  $\sigma(H)$ which are not in 
$\sigma_p(H)$ or $\sigma_c(H)$. The spectral radius of $H$ is
\beq\label{spr}
r_\sigma(H)\defeq\sup_{\lambda\in \sigma(H)} |\lambda|.
\eeq
The operator $H$ has a matrix representation $\textbf{H}$ in the
standard basis $\textbf{e}_{ik}\defeq \delta_{ik}$, where $\delta_{ik}$ is the 
Kronecker symbol.

We shall also use two standard results: first, if the operator $H$ is
bounded or closed and has a matrix representation $\textbf{H}$, then the transpose matrix $\textbf{H}^t$ is the matrix
representation of the operator $H^t:l^\infty\to l^\infty$ and (see
e.g.\ \cite{TaLa80}, \cite[Corollary II.5.3]{Gol66}
or \cite[Theorems 3.2 and 3.3]{Hal55})
\beq\label{thms}
\sigma_p(H^t)\subseteq \sigma_p(H)\cup
\sigma_r(H),
\quad \sigma_r(H)\subseteq \sigma_p(H^t).
\eeq
Second, if $H$ is bounded, then (see for example \cite[(3-5)]{TaLa80})
\beq\label{spr1}
r_\sigma(H)=\lim_{k\to\infty}\|H^k\|_1^{1/k}.
\eeq

Our aim is to classify spectra of two classes of generalized Fibonacci
and Fibonacci-like operators. For the first class of operators their spectral radii are expressed in terms
of largest real positive roots of certain polynomials and the coefficients of their
powers behave like generalized Fibonacci sequences, as we
shall see in section 2. 

The second class of operators, which also has applications in
mathematical biology, is analyzed in a similar manner in section 3.

\section{Generalized Fibonacci operators}
\label{genfib}

Let the linear operator $F_n:l^1\to l^1$ be defined by
\beq\label{rhodef}
(x_1,x_2,x_3,\cdots)\to \left(\sum_{k=n+1}^{\infty}  x_k, x_1, x_2,
  x_3, \cdots \right), \qquad n=1,2,3,\ldots
\eeq
Each $F_n$ is bounded and its matrix representation in the standard
basis is
\beq
\textbf{F}_n 
= 
\begin{array}{l}
\quad \overbrace{\phantom{xxxxxxxxxx}}^n \\
\left( \begin{array}{cccccccccc}
0 &     0 &     \cdots  &0   &1       &1      &1  & 1  & 1 & \cdots      \\
 1 &0    & 0      & 0   & 0      & 0     & 0      & 0   & 0     & \cdots\\
 0 & 1   & 0      & 0   & 0      & 0     & 0      & 0   & 0     & \cdots\\
 \vdots & \vdots   & \ddots      & \vdots   & \vdots      & \vdots
 & \vdots     & \vdots   & \vdots     & \cdots\\
 0 & 0   & 0      & 1   & 0      & 0     & 0      & 0   & 0     & \cdots\\
 0 & 0   & 0      & 0   & 1      & 0     & 0      & 0   & 0     & \cdots\\
 \vdots & \vdots         & \vdots         & \vdots      & \vdots         & \ddots        & \vdots         & \vdots      & \vdots        & \ddots\\
\end{array}\right), \label{Hdef}
\end{array}
\eeq

Following the analysis of the spectrum of $F_1$ by Halberg \cite{Hal55},
the spectrum of $F_n$ is classified in several steps which are summarized as
follows:
\begin{enumerate}
\item first, by solving the equation 
\beq\label{e1}
(\lambda I-F_n)\, x=0, \quad x\neq 0,
\eeq
we show that the point spectrum is 
\beq\label{e1a}
\sigma_p(F_n)=\{ \lambda\in\mathbb{C}: \
\lambda^{n+1}-\lambda^{n}-1=0, \ |\lambda|>1 \},
\eeq
\item second, by solving the equation 
\beq\label{e2}
(\lambda I-F_n)\, x=y, \quad x\neq 0,
\eeq
we compute the inverse $(\lambda I-F_n)^{-1}$ and show that the
resolvent set consists of all $\lambda$ such that
$|\lambda|>1$ which are not in $\sigma_p(F_n)$, that is, 
\beq\label{e3}
\rho(F_n)=\{ \lambda\in \mathbb{C}: \
 |\lambda|>1, \ \lambda^{n+1}-\lambda^{n}-1\neq 0 \},
\eeq
\item third, we analyze the transposed operator $F_n^t$ and show that 
\beq\label{e4}
\sigma_p(F_n^t)= \{ \lambda\in \mathbb{C}: |\lambda|\leq 1, \
\lambda\neq 1\},
\eeq
which, together with (\ref{thms}), implies that the residual spectrum
of $F_n$ is
\beq\label{e5}
\sigma_r(F_n)=\{ \lambda\in \mathbb{C}: \
 |\lambda|\leq 1, \ \lambda\neq 1 \}.
\eeq
\item Finally, since the spectrum of $F_n$ is closed, is also contains
  the point $\lambda =1$. Since this point is neither in the point
  spectrum nor in the residual spectrum, it must be in the continuous
  spectrum, that is 
\beq\label{e5a}
\sigma_c(F_n)=\{  1 \}.
\eeq
\end{enumerate}

We proceed with the detailed analysis of each step.

{\em Step 1.}
The equation (\ref{e1}) can be written as
\begin{align}\label{e6}
0&=\lambda x_1-x_{n+1}-x_{n+2}-x_{n+3}-\cdots, \nn \\
x_1&=\lambda x_2,\nn \\
x_2&=\lambda x_3,\nn \\
&\ \, \vdots  \\
x_k&=\lambda x_{k+1}, \nn\\
&\ \, \vdots \nn
\end{align}
Since $\lambda =0$ implies $x=0$, zero is not an element of 
$\sigma_p(F_n)$.
If $\lambda\neq 0$, by applying (\ref{e6}) recursively, we have
\beq\label{x1}
x_{k+1}=\frac{1}{\lambda}\, x_{k}=\frac{1}{\lambda^2}\, x_{k-1} =
\frac{1}{\lambda^3}\,x_{k-2}=\cdots =\frac{1}{\lambda^k}\, x_1,\qquad
k\geq 1.
\eeq
Thus
\beq\label{evec}
x=x_1 \begin{pmatrix} 1& \frac{1}{\lambda} & \frac{1}{\lambda^2} & \cdots &
  \frac{1}{\lambda^k} & \cdots
\end{pmatrix}^t
\eeq
and
\beq\label{x1a}
\| x\|_1=|x_1| \sum \frac{1}{|\lambda|^k}.
\eeq
If $|\lambda|\leq 1$, then $\| x\|_1=\infty$, so $x\notin l^1$. If
$|\lambda|>1$, then $\| x\|_1 = |x_1|\, |\lambda|/(|\lambda|-1)$.
Inserting (\ref{x1}) into the first equality of (\ref{e6}) gives
\begin{align*}
0& =\lambda x_1-x_{n+1}-x_{n+2}-x_{n+3}-\cdots \\ & =\lambda \,x_1
-\frac{1}{\lambda^n}\, x_1 -\frac{1}{\lambda^{n+1}}\, x_1 - \frac{1}{\lambda^{n+2}}\,
  x_1- \cdots \\
&= x_1 \bigg[\lambda -\frac{1}{\lambda^n} \bigg( 1
    +\frac{1}{\lambda}+\frac{1}{\lambda^2}+
\frac{1}{\lambda^3}+\cdots \bigg) \bigg]\\
&= x_1\bigg(\lambda-\frac{1}{\lambda^n}\frac{1}{1-\frac{1}{\lambda}}\bigg)\\
&=x_1 \frac{\lambda^{n+1}-\lambda^{n}-1}{\lambda^{n-1}(\lambda -1)}.
\end{align*}
Since $x_1\neq 0$, we conclude that $\sigma_p(F_n)$ consists of those
roots of the polynomial
\beq\label{poly}
p_{n+1}(\lambda)\defeq \lambda^{n+1}-\lambda^{n}-1
\eeq
for which $|\lambda|>1$, as stated in (\ref{e1a}).\footnote{These
  roots are the eigenvalues and the vectors $x$ defined by 
  (\ref{evec}) are the corresponding eigenvectors.}

\begin{figure}[htb]
  \centering
\begin{tabular}{cc}
  \includegraphics[scale=0.7]{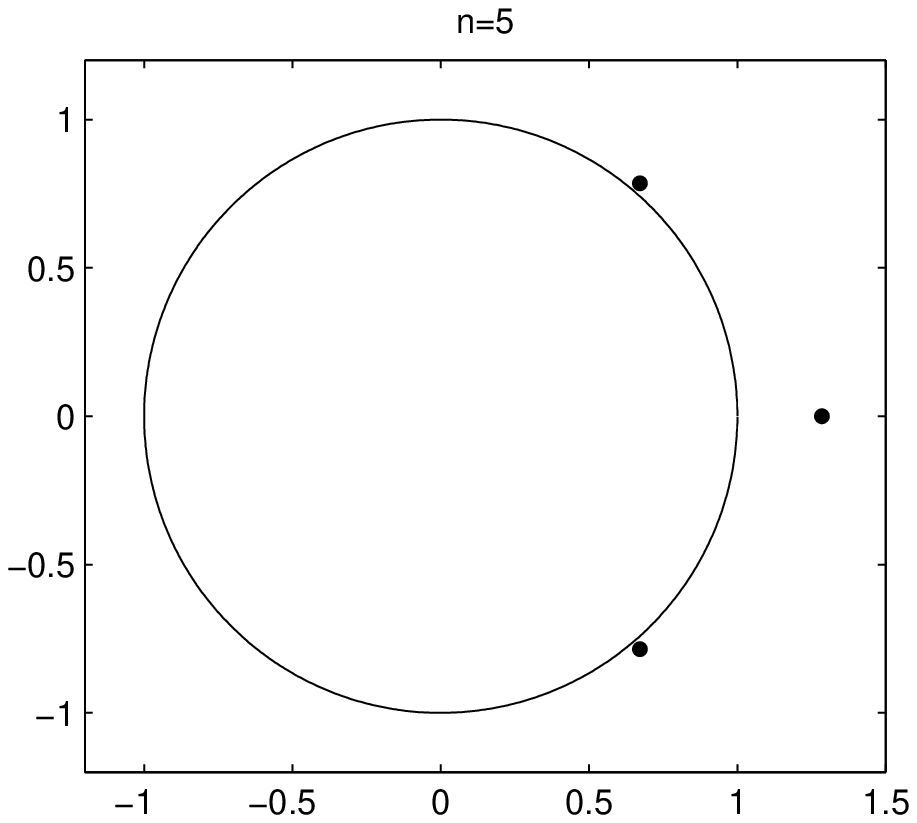} &
  \includegraphics[scale=0.7]{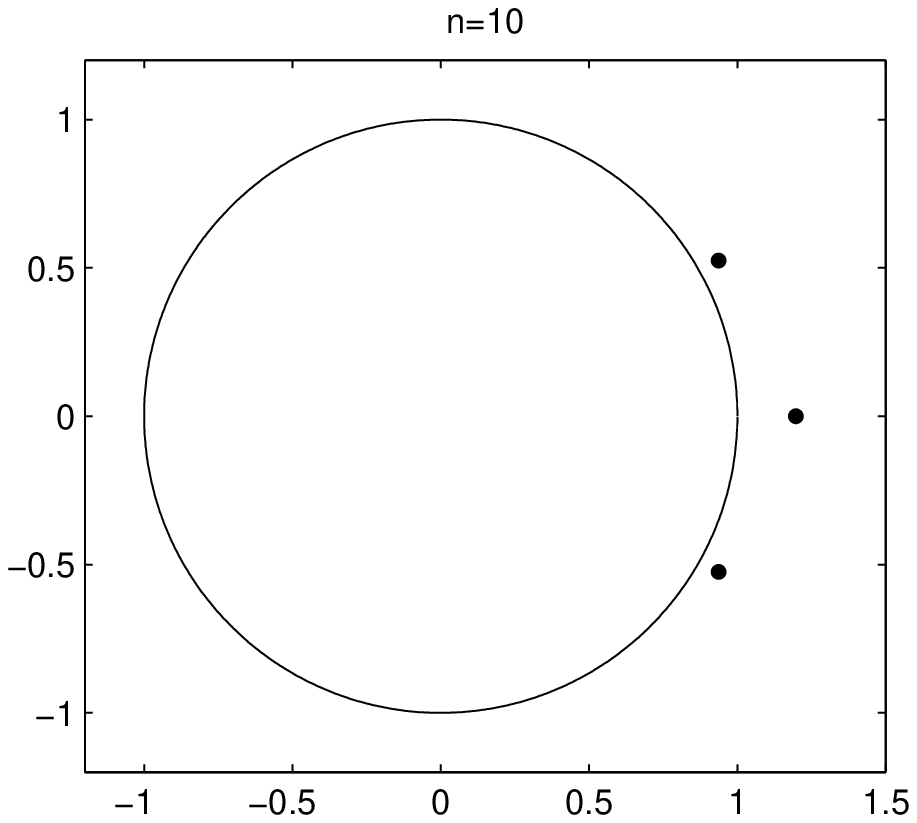}\\
  \includegraphics[scale=0.7]{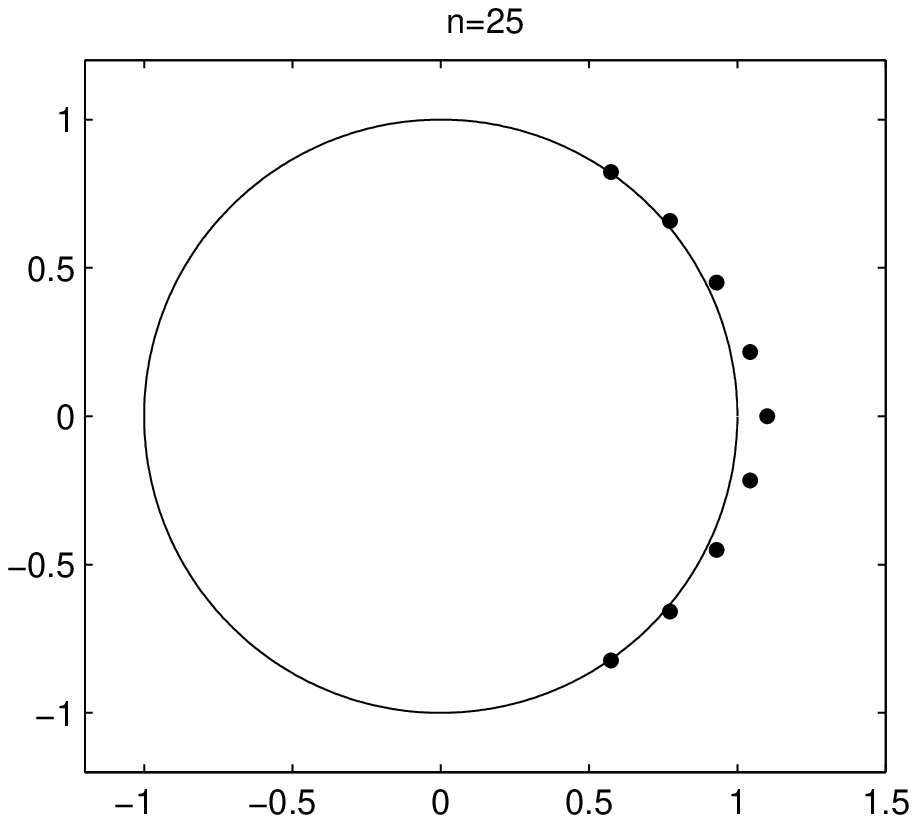} &
  \includegraphics[scale=0.7]{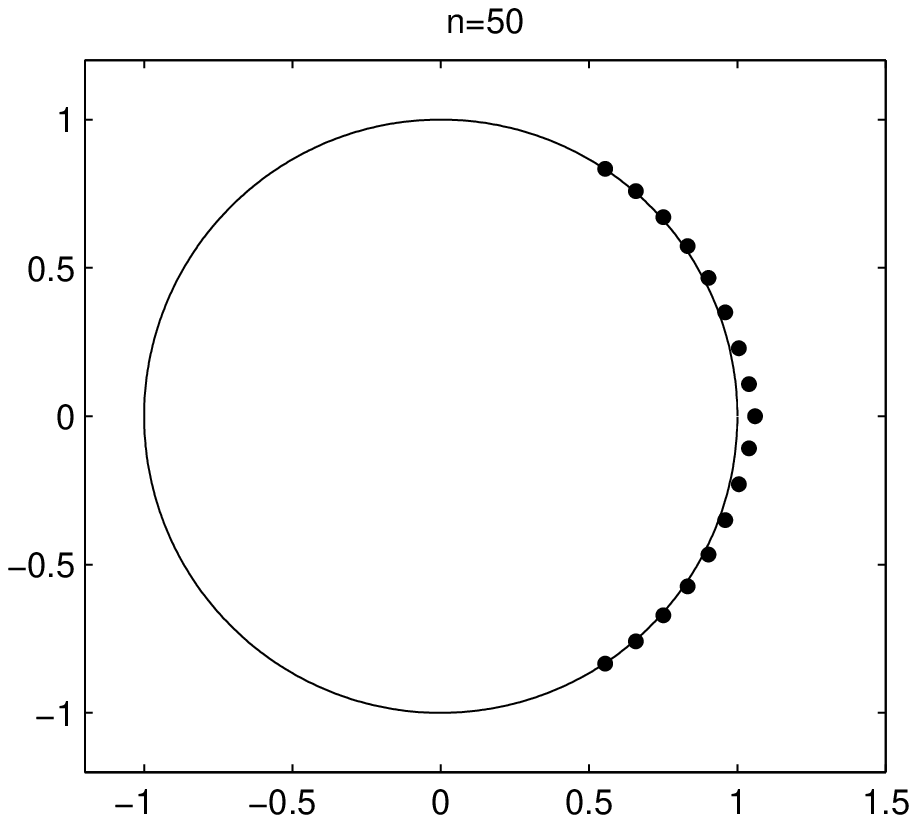}
\end{tabular}
\caption{The point spectra $\sigma_p(F_n)$ for various values of $n$. }
  \label{fig:1}
\end{figure}

Since $p_{n+1}(1)=-1<0$ and $p'_{n+1}(\lambda)>0$ for 
$\lambda\in\mathbb{R}, \lambda\geq 1$,
that is, $p_{n+1}$ is strictly increasing for $\lambda>1$,
we conclude that $F_n$ has exactly one real eigenvalue larger
than one. Let us denote this eigenvalue by $\la_{\max}(F_n)$.
By Ostrovsky's theorem \cite[Theorem 1.1.4, p. 3]{Pra04},
$\la_{\max}(F_n)$ is the unique positive root of $p_{n+1}(\lambda)$
and the absolute values of all other roots are strictly smaller.
Consequently, all other eigenvalues of $F_n$ are in absolute value
strictly smaller than $\la_{\max}(F_n)$ which, in turn, implies 
\beq\label{rla}
r_\sigma(F_n)=\la_{\max}(F_n).
\eeq
Figure \ref{fig:1} shows $\sigma_p(F_n)$ for various values of $n$.

{\em Step 2.}
The equation (\ref{e2}) can be written as
\begin{align}\label{s2}
y_1& =\la x_1 -x_{n+1}-x_{n+2}-x_{n+3}-\cdots, \nn \\
x_2&=\frac{1}{\la}\, (x_1+y_2),\nn \\
x_3&=\frac{1}{\la}\, (x_2+y_3),\\
&\ \, \vdots  \nn \\
x_{k+1}&=\frac{1}{\la}\, (x_k+y_{k+1}),\nn \\
&\ \, \vdots \nn
\end{align}
By setting 
\[
u=\sum_{k=n+1}^{\infty} x_k, \qquad 
v=\sum_{k=n+1}^{\infty} y_k,
\]
and using (\ref{s2}), we have
\[
u=\frac{1}{\la}\, x_n + \frac{1}{\la}\, u + \frac{1}{\la}\, v
, \qquad y_1=\la x_1-u.
\]
After rearranging, we have
\[
u=\frac{1}{\la -1}(x_n+v).
\]
Thus, 
\beq\label{u1}
y_1=\la x_1 - \frac{1}{\la -1}(x_n+v).
\eeq
By recursively applying (\ref{s2}), we have
\begin{align}\label{u2}
x_2&=\frac{1}{\la}\, (x_1+y_2),\nn \\
x_3&=\frac{1}{\la}\, (x_2+y_3)=\frac{1}{\la^2}\, x_1+ \frac{1}{\la^2}\, y_2 
+ \frac{1}{\la}\,y_3, \nn \\
&\ \, \vdots  \nn \\
x_{k+1}&=\frac{1}{\la}\, (x_k+y_{k+1})=\frac{1}{\la^k}\, x_1+\frac{1}{\la^k}\,
y_2+\frac{1}{\la^{k-1}}\, y_3+\frac{1}{\la^{k-2}}\, y_4+\cdots +
\frac{1}{\la}\, y_{k+1},\\
&\ \, \vdots \nn
\end{align}
Inserting $x_n$ into (\ref{u1}) gives
\[
y_1=\la \, x_1 -\frac{1}{\la-1}\bigg(
\frac{1}{\la^{n-1}}\, x_1  +\frac{1}{\la^{n-1}}\, y_2+
 \frac{1}{\la^{n-2}}\, y_3 + \cdots + \frac{1}{\la}\, y_n+v\bigg),
\]
and solving for $x_1$ gives
\[
x_1=\frac{1}{\la^{n+1}-\la^{n}-1}
\bigg(  \la^{n-1}(\la -1)\,  y_1 + y_2 + \la\, y_3 + \la^2\,  y_4+
\cdots + \la^{n-2}\, y_n + \la^{n-1}\, v\bigg).
\] 
By inserting this into (\ref{u2}) we have
\[
x=(\lambda I-F_n)^{-1}y=\frac{1}{\la^{n+1}-\la^{n}-1} \, (A+B)\,y
\]
where the matrix representations of $A$ and $B$ are
given by\footnote{Next row of $\textbf{A}$ is obtained by dividing the
  previous row by $\lambda$.}
\[
\textbf{A}
=
\left( \begin{array}{cccccccccc}
(\la-1)\la^{n-1} &     1  &     \la  &\la^2   &\la^3 &\cdots       &\la^{n-2}
&\la^{n-1}  & \la^{n-1}  & \cdots      \\
(\la-1)\la^{n-2} &  \frac{1}{\la}  &    1  &\la   &\la^2  &\cdots
&\la^{n-3}  & \la^{n-2}  & \la^{n-2} & \cdots      \\
 \vdots & \vdots   & \vdots      & \vdots   & \vdots      & \vdots
 & \vdots     & \vdots   & \vdots     & \cdots\\
(\la-1)\la &  \frac{1}{\la^{n-2}}  &  \frac{1}{\la^{n-3}}  &
\frac{1}{\la^{n-4}} & \cdots & \frac{1}{\la}   &1        &\la
&\la & \cdots      \\
(\la-1) &  \frac{1}{\la^{n-1}}  &  \frac{1}{\la^{n-2}}  &
\frac{1}{\la^{n-3}} & \cdots & \frac{1}{\la^2}   &\frac{1}{\la}        &1
&1 & \cdots      \\
(\la-1)\,\frac{1}{\la} &  \frac{1}{\la^{n}}  &  \frac{1}{\la^{n-1}}  &
\frac{1}{\la^{n-2}} & \cdots & \frac{1}{\la^3}   &\frac{1}{\la^2}        &\frac{1}{\la}
&\frac{1}{\la} & \cdots      \\
 \vdots & \vdots         & \vdots         & \vdots      & \vdots         & \ddots        & \vdots         & \vdots      & \vdots        & \ddots\\
\end{array}\right),
\]
and 
\[
\textbf{B}=
\left( \begin{array}{cccccc}
0 & 0 & 0 & 0 & 0 & \cdots \\
0 & \frac{1}{\la} & 0 & 0 & 0 & \cdots \\
0 & \frac{1}{\la^2} & \frac{1}{\la} & 0 & 0 & \cdots \\
0 & \frac{1}{\la^3} & \frac{1}{\la^2} & \frac{1}{\la} & 0 & \cdots \\
\vdots & \vdots & \vdots & \vdots & \vdots & \ddots
\end{array}\right),
\]
respectively. Obviously, for $|\la|>1$ we have $\|A\|_1< \infty$ and 
$\|B\|_1< \infty$. Thus, for $|\la|>1$ and $\la$ not being the root of 
$\la^{n+1}-\la^n-1$, the operator $\la I-F_n$ has a bounded inverse, so
the resolvent set of $F_n$ is given by 
(\ref{e3}).

{\em Step 3.}
The point spectrum of the transposed operator $F_n^t$
consists of all $\lambda\in\mathbb{R}$ such that
\[
(\lambda I - F_n^t)\,x=0, \qquad x\neq 0,\quad \|x\|_{\infty}< \infty.
\]
This is equivalent to
\begin{align*}
x_2&=\la\, x_1, \\
x_3&=\la \, x_2 =\la^2\, x_1, \\
&\ \, \vdots  \\
x_{n+1}&=\la \, x_n =\la^{n}\, x_1, \\
x_{n+2}&=\la \, x_{n+1} -x_1=(\la^{n+1}-1)\, x_1, \\
x_{n+3}&=\la \, x_{n+2} -x_1=(\la^{n+2}-\la-1)\, x_1, \\
&\ \, \vdots  \\
x_{k}&=(\la^{k-1}-\la^{k-n-2}- \la^{k-n-3}-\cdots - \la-1)\, x_1, \\
&\ \, \vdots
\end{align*}
Therefore,
\[
x_k = \bigg(\la^{k-1}-\frac{\la^{k-n-1}-1}{\la -1}\bigg)\, x_1. \\
\]
For $|\la|\leq 1$, $\la\neq 1$ we have
\[
|x_k|< \bigg(1+\frac{2}{|\la-1|}\bigg) |x_1|,
\]
which implies $\|x\|_\infty < \infty$. For $\la =1$ we have
\begin{align*}
x_2&=x_1,\\
x_3&=x_1,\\
&\ \, \vdots  \\
x_{n+1}&=x_1,\\
x_{n+2}&=0, \\
x_{n+3}&=-x_1, \\
x_{n+4}&=-2\, x_1, \\
&\ \, \vdots  \\
x_{k}&=-(k-n-2)\, x_1, \\
&\ \, \vdots
\end{align*}
so $\|x\|_\infty = \infty$. We conclude that the point spectrum of
$F_n^t$ is given by (\ref{e4}). This, in turn, implies (\ref{e5}) and
(\ref{e5a}) as described before.

\subsection{Relationship to generalized Fibonacci sequences}

In this section we describe the relationship between operators $F_n$
and generalized Fibonacci sequences. 
A generalized Fibonacci sequence $\{ f^{(n)}\}$ is defined by
\beq\label{d1}
f^{(n)}_1=1, \quad f^{(n)}_2=1, \cdots,  f^{(n)}_{n+1}=1,\quad
f^{(n)}_k=f^{(n)}_{k-1}+f^{(n)}_{k-n-1}, \qquad k > n+1.
\eeq
For $n=1$ this definition yields the classical Fibonacci sequence
\beq\label{fib}
f_1=1, \quad f_2=1,\quad f_k=f_{k-1}+f_{k-2}, \quad k>2.
\eeq
By induction we can prove that the $k$-th power of the matrix
$\textbf{F}_n$ from (\ref{Hdef}) for $k>n$ has the form
\[
\textbf{F}_n^k=
 \left( \begin{array}{ccccccccc}
f^{(n)}_{k-n} &     f^{(n)}_{k-n+1} &     f^{(n)}_{k-n+2}  &\cdots
&f^{(n)}_{k-1} & f^{(n)}_{k}     & f^{(n)}_{k}  & f^{(n)}_{k} &\cdots
\\
f^{(n)}_{k-n-1} &     f^{(n)}_{k-n} &     f^{(n)}_{k-n+1}  &\cdots
&f^{(n)}_{k-2} & f^{(n)}_{k-1}     & f^{(n)}_{k-1}  & f^{(n)}_{k-1} &\cdots
\\
 \vdots & \vdots & \vdots         & \vdots      & \vdots
 & \vdots        & \vdots         & \vdots  & \cdots \\
f^{(n)}_{1} &     f^{(n)}_{2} &     f^{(n)}_{3}  &\cdots
&f^{(n)}_{n} & f^{(n)}_{n+1}     & f^{(n)}_{n+1}  & f^{(n)}_{n+1} &\cdots
\\
0 & 1    & 1    & \cdots  & 1      & 1     & 1      & 1    & \cdots\\
0 & 0    & 1    & \cdots  & 1      & 1     & 1      & 1    & \cdots\\
\vdots & \vdots  & \vdots & \vdots & \vdots & \vdots  & \vdots & \vdots &\cdots   \\
0 & 0    & 0    & \cdots  & 0      & 1     & 1      & 1    & \cdots \\
1 & 0    & 0    & \cdots  & 0      & 0     & 0      & 0    & \cdots \\
0 & 1    & 0    & \cdots  & 0      & 0     & 0      & 0    & \cdots \\
0 & 0    & 1    & \cdots  & 0      & 0     & 0      & 0    & \cdots \\
\vdots & \vdots  & \vdots & \ddots & \vdots & \vdots  & \vdots & \vdots &\cdots   \\
\end{array}\right).
\]
We conclude that
\beq\label{nr1}
\| F_n^k\|_1=1+\sum_{i=1}^k f^{(n)}_i.
\eeq
By applying (\ref{d1}) to the terms in parentheses we have
\begin{align*}
2\,\sum_{i=1}^k f^{(n)}_i & =  f^{(n)}_1+\cdots + f^{(n)}_n+
\big( f^{(n)}_{n+1}+  f^{(n)}_{n+2}+\cdots+   f^{(n)}_{k-1}\big)+
 f^{(n)}_k \\
&\quad + \big( f^{(n)}_1+ f^{(n)}_2+\cdots +  f^{(n)}_{k-n-1}\big)
+  f^{(n)}_{k-n}+ f^{(n)}_{k-n+1}+\cdots +  f^{(n)}_k\\
&= f^{(n)}_1+f^{(n)}_2+\cdots +f^{(n)}_n+f^{(n)}_{n+2}+ \cdots
+f^{(n)}_{k}\\ & \quad +
f^{(n)}_k+f^{(n)}_{k-n}+f^{(n)}_{k-n+1}+\cdots+f^{(n)}_k\\
&= \sum_{i=1}^k f^{(n)}_i - f^{(n)}_{n+1} 
+f^{(n)}_{k-n}+f^{(n)}_{k-n+1}+\cdots+f^{(n)}_k +f^{(n)}_k.
\end{align*}
From this, by applying  (\ref{d1}) again recursively, we obtain
\begin{align*}
\sum_{i=1}^k f^{(n)}_i&=f^{(n)}_{k-n}+f^{(n)}_{k-n+1}+\cdots+f^{(n)}_k
+f^{(n)}_k-1 \\
&= f^{(n)}_{k-n+1}+\cdots+f^{(n)}_k
+f^{(n)}_{k+1}-1 \\ &=  f^{(n)}_{k-n+2}+\cdots+f^{(n)}_{k+1}
+f^{(n)}_{k+2}-1 \\
&\ \, \vdots \\ 
&= f^{(n)}_{k+n+1}-1.
\end{align*}
Inserting this into (\ref{nr1}) gives
\beq\label{nr2}
\| F_n^k\|_1= f^{(n)}_{k+n+1}
\eeq
and from (\ref{rla}) it follows that 
\[
\lim_{k\to\infty}  \big(f^{(n)}_{k+n+1}\big)^{1/k}=\la_{\max}(F_n).
\]
Also, by using standard techniques in analyzing linear recurrence
relations with constant coefficients, we can prove that for all $i,j$
\footnote{The proof is derived using the fact that $f^{(n)}_l$ has the
                 form 
  $f^{(n)}_l=\alpha\,\la^l_{\max}(F_n) +\sum_{i=1}^{n} \alpha_i \la_i^l$, where
  $\la_{\max}(F_n)$ and $\la_i$ are
  the roots of the characteristic polynomial (\ref{poly}), and 
  $|\la_{\max}(F_n)|>|\la_i|$.}
\[
\lim_{k\to\infty}
\frac{[\textbf{F}_n^k]_{i,j}}{[\textbf{F}_n^k]_{i+1,j}}
\equiv 
\lim_{m\to\infty} \frac{f^{(n)}_{m+1}}{f^{(n)}_{m}}=\la_{\max}(F_n).
\]

For example, by setting $n=1$ we have for the Fibonacci sequence (\ref{fib})
\begin{align*}
\lim_{k\to\infty}
(f_{k+2})^{1/k}& =r_\sigma(F_1)=\frac{1+\sqrt{5}}{2}, \\
\lim_{k\to\infty} \frac{f_{k+1}}{f_{k}}&=\frac{1+\sqrt{5}}{2}.
\end{align*}

\section{Fibonacci-like operators}

Now we would like to consider the family of linear operators  
$\Gamma_n:l^1\to l^1$ defined by
\beq\label{rhodefg}
(x_1,x_2,x_3,\cdots)\to \left(\rho \sum_{k=n+1}^{\infty}  (k-n)\, x_k, x_1, x_2,
  x_3, \cdots \right), \qquad n=1,2,3,\ldots
\eeq
for some real positive $\rho$.
The domain of $\Gamma_n$ is 
\[
\textrm{Dom } \Gamma_n = \left\{ x\in l^1 : \left|\sum_{k=n+1}^{\infty} 
(k-n) x_k \right|  <\infty \right\},
\]
and its matrix representation in the standard basis is
\beq
\mathbf{\Gamma}_n 
= 
\begin{array}{l}
\quad\, \overbrace{\phantom{xxxxxxxxxx}}^n \\
\left( \begin{array}{cccccccccc}
0 &     0 &     \cdots  &0   &\rho       &2\rho      &3\rho  & 4\rho
& 5\rho  & \cdots      \\
 1 &0    & 0      & 0   & 0      & 0     & 0      & 0   & 0     & \cdots\\
 0 & 1   & 0      & 0   & 0      & 0     & 0      & 0   & 0     & \cdots\\
 \vdots & \vdots   & \ddots      & \vdots   & \vdots      & \vdots
 & \vdots     & \vdots   & \vdots     & \cdots\\
 0 & 0   & 0      & 1   & 0      & 0     & 0      & 0   & 0     & \cdots\\
 0 & 0   & 0      & 0   & 1      & 0     & 0      & 0   & 0     & \cdots\\
 \vdots & \vdots         & \vdots         & \vdots      & \vdots         & \ddots        & \vdots         & \vdots      & \vdots        & \ddots\\
\end{array}\right). \label{HdefGa}
\end{array}
\eeq

However, the operator $\Gamma_n$ is not 
closed as illustrated by the following example.

\begin{example}\label{ex1}
Let us define the sequence $\{x^{(m)}\}$ of vectors in $l^1$ by
\[
\begin{array}{l}
\qquad \qquad  \overbrace{\phantom{xxxxxxx}}^{m+n-1} \\
x^{(m)}= \begin{pmatrix}0 &\cdots
  &0&\displaystyle\frac{1}{m}&0&\cdots \end{pmatrix}^t.
\end{array}
\]
Then
\[
x^{(m)}
 \to \begin{pmatrix}0 & 0 &\cdots
 \end{pmatrix}^t,
 \]
while
\[
\mathbf{\Gamma}_n \, x^{(m)} = \begin{pmatrix}\rho&0&\cdots&0&\frac{1}{m}&0&\cdots
\end{pmatrix}^t
\to \begin{pmatrix} \rho &0,&0&\cdots
\end{pmatrix}^t.
\]
\end{example}

Although the point spectrum of $\Gamma_n$ is defined and can be computed
in a standard manner (see later), the resolvent set of $\Gamma_n$ is empty,
which makes the analysis of $\Gamma_n$ less interesting. Instead, we shall
consider the family of operators $G_n:l_1\to l_1$ formally defined by
\[
G_n=D_n \Gamma_n D_n^{-1},
\]
where
\[
D_n=\diag \big(\overbrace{1,\cdots,1}^n,\rho,2\rho, 3\rho, 4\rho,\cdots \big).
\]
That is, for $n\in\mathbb{N}$ the operator $G_n$ is defined by
\[
(x_1,x_2,x_3,\cdots)\to \left(\sum_{k=n+1}^{\infty} x_k, x_1,
  x_2,\cdots,x_{n-1}, \rho\,x_n,2\, x_{n+1},\frac{3}{2}\,x_{n+2},
 \frac{4}{3}\,x_{n+3}, \frac{5}{4}\,x_{n+4},
 \cdots \right),
\]
and its matrix representation in the standard basis is 
\beq
\textbf{G}_n 
= 
\begin{array}{l}
\quad\, \overbrace{\phantom{xxxxxxxxxx}}^n \\
\left( \begin{array}{cccccccccc}
0 &     0 & \cdots  & 0   & 1      & 1     & 1    & 1   & 1     & \cdots\\
1 &     0 & \cdots  & 0   & 0      & 0     & 0    & 0   & 0     & \cdots\\
0 &     1 & \cdots  & 0   & 0      & 0     & 0    & 0   & 0     & \cdots\\
\vdots & \vdots   & \ddots   & \vdots   & \vdots      & \vdots
 & \vdots     & \vdots   & \vdots     & \cdots\\
0 &     0 & \cdots  &\rho & 0      & 0     & 0    & 0   & 0     & \cdots\\
0 &     0 & \cdots  & 0   & 2      & 0     & 0    & 0   & 0     & \cdots\\
0 &     0 & \cdots  & 0   & 0      & \frac{3}{2}    & 0   & 0   & 0     & \cdots\\
0 &     0 & \cdots  & 0   & 0      & 0   & \frac{4}{3}   & 0   & 0     &
\cdots\\
0 &     0 & \cdots  & 0   & 0      & 0   & 0  & \frac{5}{4}   & 0     &
\cdots\\
  \vdots & \vdots         & \vdots         & \vdots      & \vdots
  & \vdots         & \vdots      & \vdots        & \ddots & \vdots
\end{array}\right). \label{Hdefg}
\end{array}
\eeq


Let us define the polynomial $q_{n+1}(\la)$ by
\beq\label{polq}
q_{n+1}(\lambda)=
  \lambda^{n+1}-2\,\lambda^{n}+\lambda^{n-1}-\rho.
\eeq
Similarly as in section \ref{genfib}, the spectrum of $G_n$ is classified in several steps which are summarized as
follows:
\begin{enumerate}
\item first, by solving the equation 
\beq\label{e1g}
(\lambda I-G_n)\, x=0, \quad x\neq 0,
\eeq
we show that the point spectrum is 
\beq\label{e1ag}
\sigma_p(G_n)=\{ \lambda\in\mathbb{C}: \ q_{n+1}(\la)=0, \ |\lambda|>1
\},\quad n\geq 2.
\eeq

\item second, by solving the equation 
\beq\label{e2g}
(\lambda I-G_n)\, x=y, \quad x\neq 0,
\eeq
we compute the inverse $(\lambda I-G_n)^{-1}$ and show that the
resolvent set consists of all $\lambda$ such that
$|\lambda|>1$ which are not in $\sigma_p(G_n)$,
\beq\label{e3g}
\rho(G_n)=\{ \lambda\in \mathbb{C}: \
 |\lambda|>1, \ \lambda\notin\sigma_p(G_n) \},
\eeq



\item third, we analyze the transposed operator $G_n^t$ and show that 
\beq\label{e4G}
\sigma_p(G_n^t)= \{ \lambda\in \mathbb{C}: |\lambda|\leq 1, \
\lambda\neq 1\},
\eeq
which, together with (\ref{thms}), implies that the residual spectrum
of $G_n$ is
\beq\label{e5G}
\sigma_r(G_n)=\{ \lambda\in \mathbb{C}: \
 |\lambda|\leq 1, \ \lambda\neq 1 \}.
\eeq
\item Finally, since the spectrum of $G_n$ is closed, is also contains
  the point $\lambda =1$. Since this point is neither in the point
  spectrum nor in the residual spectrum, it must be in the continuous
  spectrum, that is 
\beq\label{e5aG}
\sigma_c(G_n)=\{  1 \}.
\eeq
\end{enumerate}
The proofs are similar to the ones from section \ref{genfib}, but more
tedious.

{\em Step 1.}
The equation (\ref{e1g}) can be written as
\begin{align}\label{e6g}
0&=\lambda x_1-x_{n+1}-x_{n+2}-x_{n+3}-\cdots,  \\
x_k&=\lambda x_{k+1}, \qquad k=1,\cdots,n-1, \nn\\
\rho \, x_n&= \lambda\, x_{n+1}, \nn \\
\frac{k-n+1}{k-n}\,x_k&=\lambda\, x_{k+1}, \qquad k=n+1,n+2,\cdots.\nn
\end{align}
Since $\lambda =0$ implies $x=0$, zero is not an element of 
$\sigma_p(G_n)$.
If $\lambda\neq 0$, by applying (\ref{e6g}) recursively, we obtain
\begin{align}\label{x1A}
x_k&=\frac{1}{\lambda^{k-1}}\, x_1, \qquad k=2,3,\cdots,n,\nn \\
x_k&=\rho\, \frac{k-n}{\lambda^{k-1}}\, x_1,\qquad k=n+1,n+2,\cdots.
\end{align}
This, in turn, implies
\[
x=x_1 \big(\begin{array}{ccccccccccc}1&\frac{1}{\lambda}& \frac{1}{\lambda^2}&
\cdots &  \frac{1}{\lambda^{n-1}}& \frac{\rho}{\lambda^n}&
 \frac{2\,\rho}{\lambda^{n+1}}& \frac{3\,\rho}{\lambda^{n+2}}&
\cdots &  \frac{(k-n)\,\rho}{\lambda^{k-1}} & \cdots
\end{array}\big)^t,
\]
so that 
\[
\|x\|_1 =|x_1| \,
\bigg(\frac{1-(\frac{1}{\lambda})^n}{1-\frac{1}{\lambda}}+\frac{\rho}{\lambda^n}
(1+\frac{2}{\lambda}+\frac{3}{\lambda^2}+\cdots )\bigg). 
\]
%
If $|\lambda|\leq 1$, then $\| x\|_1=\infty$, so $x\notin l^1$. If
$|\lambda|>1$, then, by using differentiation of the geometric series, we
have
\[
\| x\|_1 = |x_1| \left(
\frac{\lambda^n-1}{\lambda^{n-1}(\lambda-1)}
+\frac{\rho}{\lambda^n}\, \frac{\lambda^2}{(\lambda-1)^2}\right)<\infty,
\]
thus, $x\in l^1$.
By inserting (\ref{x1A}) into the first equality of (\ref{e6g}) and
using differentiation of the geometric series we have
\begin{align*}
0& 
=\lambda \,x_1
-\rho\, \bigg( \frac{1}{\lambda^n}\, x_1 +\frac{2}{\lambda^{n+1}}\, x_1 + \frac{3}{\lambda^{n+2}}\,
  x_1+ \cdots \bigg)\\
&= x_1 \bigg[\lambda -\frac{\rho}{\lambda^n} \bigg( 1
    +\frac{2}{\lambda}+\frac{3}{\lambda^2}+
\frac{4}{\lambda^3}+\cdots \bigg) \bigg]\\
&=
x_1\bigg(\lambda-\frac{\rho}{\lambda^n}\frac{1}{(1-\frac{1}{\lambda})^2}\bigg).
\end{align*}
Finally, solving this equation with $x_1\neq 0$ gives (\ref{e1ag}).

We shall now prove that $\sigma_p(G_n)$ consists of $\la_{\max}(G_n)$,
a unique simple real eigenvalue larger than one and all other
eigenvalues have modulus smaller than  $\la_{\max}(G_n)$. This also implies
\beq\label{rlag}
r_\sigma(G_n)=\la_{\max}(G_n).
\eeq
The proof is based on the ideas from the proof of \cite[Theorem 1.1.4, pp. 3-4]{Pra04}.
Indeed, if $n=1$ then the roots of $q_2(\lambda)$ are $1\pm
\sqrt{\rho}$ and the statement holds. 
For $n\geq 2$ we have
\begin{align}
q_{n+1}(\lambda)&=\lambda^{n-1}(\lambda-1)^2-\rho, \label{os1}
\\
q'_{n+1}(\lambda)&=\lambda^{n-2}[(n+1)\lambda^2-2n\lambda+(n-1)].
\label{os2}
\end{align}
Since $q_{n+1}(1)=-\rho<0$ and $q'_{n+1}(\lambda)>0$ for 
$\lambda\in\mathbb{R}, \lambda> 1$,
that is, $q_{n+1}(\la)$ is strictly increasing for $\lambda>1$,
we conclude that $q_{n+1}(\lambda)$ has exactly one real root larger
than one or, equivalently, that $G_n$ has exactly one real eigenvalue larger
than one.  Let us denote this eigenvalue by $\la_{\max}(G_n)$.
Let $z\neq\la_{\max}(G_n)$ be some other real or complex eigenvalue of
$G_n$ and let $\zeta =|z|>1$. 
Since $z$ is also the root of $q_{n+1}(\lambda)$, the relation (\ref{os1}) implies
\[
z^{n-1}(z-1)^2=\rho,
\]
which, in turn, implies
\begin{equation} \label{os3}
|z|^{n-1}|z-1|^2=\rho.
\end{equation}
Since $\zeta >1$, this implies
\begin{equation*} 
\zeta^{n-1} (\zeta-1)^2\leq \rho,
\end{equation*}
or
\begin{equation}\label{os5}
q_{n+1}(\zeta)=\zeta^{n-1} (\zeta -1)^2-\rho\leq 0.
\end{equation}
Since $q_{n+1}(\lambda)$ is strictly increasing for $\lambda>1$, 
and $q_{n+1}(\la_{\max}(G_n))=0$, we conclude that
$\zeta\leq\la_{\max}(G_n)$ and that the equality in (\ref{os5}) 
holds only if
$\zeta=\la_{\max}(G_n)$. But, the equality in (\ref{os5}) and (\ref{os3}) imply
\[
|z-1|=\zeta-1,
\]
that is, $z\in\mathbb{R}$ and $z=\pm\zeta=\pm \la_{\max}(G_n)$.
The choice $z=- \la_{\max}(G_n)$ is impossible since 
$q_{n+1}(-\la_{\max}(G_n))\neq 0$, and the second choice contradicts
the assumption $z\neq\la_{\max}(G_n)$.
Therefore, $\zeta<\la_{\max}(G_n)$ as desired.

\begin{remark} Although the above analysis is sufficient for our purposes, by
inspecting the polynomial $q_{n+1}(\lambda)$ and its derivative from
(\ref{os1}) and (\ref{os2}), respectively, we can establish further
facts about its roots. 
  From (\ref{os2}) we see that the derivative $q'_{n+1}(\lambda)$ 
  has exactly two real positive simple roots
  $\lambda_1=\frac{n-1}{n+1}$ and  $\lambda_2=1$ and, if $n>2$, also the root
  $\lambda_0=0$. If $n>3$ then $\lambda_0$ is multiple.
  Let $\rho_0=4 \lambda_1^{n-1}/(n+1)^2$.
  The number of real roots of $q_{n+1}(\la)$ in the   
  open interval $(0,\la_{\max}(G_n))$ is governed by $\rho$ as follows:
  if $\rho>\rho_0$, then there are no such roots, if  $\rho=\rho_0$
  there is exactly one root equal to $\lambda_1$ and if  $\rho<\rho_0$
  there are exactly two roots, one smaller and one larger than
  $\lambda_1$.
  Finally, if $n$ is odd, then $q_{n+1}(\la)$ also has a simple
  negative real root.
  As we have already proved,
  $\lambda_{\max}(G_n)$ is the root with strictly maximal absolute value. 
\end{remark}

\begin{remark}\label{rem2}
It is easy to see that the point spectrum of $\Gamma_n$ from
(\ref{rhodefg}) and (\ref{HdefGa}) is equal to the point
spectrum of $G_n$.
\end{remark}



{\em Step 2.}
The equation (\ref{e2g}) can be written as
\begin{align}\label{s2g}
y_1& =\la x_1 -x_{n+1}-x_{n+2}-x_{n+3}-\cdots, \\
x_{k+1}&=\frac{1}{\la}\, (x_k+y_{k+1}), \qquad k=1,\cdots,n-1, \nn \\
x_{n+1}&=\frac{1}{\la}\, (\rho\,x_n+y_{n+1}), \nn \\
x_{k+1}&=\frac{1}{\la}\, \bigg(\frac{k-n+1}{k-n}\,x_k+y_{k+1}\bigg), 
\qquad k=n+1,n+2,\cdots.\nn 
\end{align}
By recursively applying the above three equalities, we obtain
\begin{align}\label{u2G}
x_k&=\frac{(k-n)\,\rho}{\la^{k-1}}\,x_1+
\frac{(k-n)\,\rho}{\la^{k-1}}\,y_2+
\frac{(k-n)\,\rho}{\la^{k-2}}\,y_3+ \cdots +
\frac{(k-n)\,\rho}{\la^{k-n+1}}\,y_n+\\
&\quad +
\frac{k-n}{1}\, \frac{1}{\la^{k-n}}\,y_{n+1}+
\frac{k-n}{2}\, \frac{1}{\la^{k-n-1}}\,y_{n+2}+
\frac{k-n}{3}\, \frac{1}{\la^{k-n-2}}\,y_{n+3}
+\cdots + \nn \\
&\quad +\frac{k-n}{k-n-1}\, \frac{1}{\la^2}\,y_{k-1}+
\frac{1}{\la}\,y_{k},\qquad k=n+1,n+2,\cdots.\nn
\end{align}
For $|\la|>1$, by inserting  (\ref{u2G}) into (\ref{s2g}), rearranging, and
using differentiation of the geometric series, we have
\begin{align}\label{yG}
y_1&=\la\,x_1-\frac{\la^2}{(\la-1)^2}\bigg(\frac{\rho}{\la^n}\,x_1+\frac{\rho}{\la^n}\,y_2+\frac{\rho}{\la^{n-1}}\,y_3 
+\cdots+
\frac{\rho}{\la^2}\,y_{n}+
\frac{1}{\la}\,y_{n+1}\bigg)- \nn
\\
&\quad -
\alpha_{n+2}\,y_{n+2}-\alpha_{n+3}\,y_{n+3}-\alpha_{n+4}\,y_{n+4}-
\cdots,
\end{align}
where
\[
\alpha_{n+k}=\frac{1}{\la}+\frac{k+1}{k}\,\frac{1}{\la^2}+
\frac{k+2}{k}\,\frac{1}{\la^3}+ 
\frac{k+3}{k}\,\frac{1}{\la^4}+\cdots,\qquad k=2,3,4,\cdots,
\]
and
\begin{equation}\label{aG}
|\alpha_{n+k}|\leq\frac{1}{|\la|}\, \bigg(1+
\sum_{i=1}^{\infty}\frac{k+i}{k}\,\frac{1}{|\la|^i}\bigg)\leq 
\frac{1}{|\la|}\, \bigg(1+
\sum_{i=1}^{\infty}(i+1)\,\frac{1}{|\la|^i}\bigg)=
\frac{|\la|}{(|\la|-1)^2}.
\end{equation}
Solving for (\ref{yG}) for $x_1$ gives
\begin{align}\label{xG}
x_1&=\frac{1}{q_{n+1}(\la)}
\bigg[ \la^{n-2}(\la-1)^2\,  y_1 + \rho\,y_2 + \rho\,\la\, y_3 + \rho\,\la^2\,  y_4+
\cdots +\rho\,\la^{n-2}\, y_{n}+\la^{n-1}\, y_{n+1}+ \nn
\\
 &\quad +\bar\alpha_{n+2}\,y_{n+2}+\bar\alpha_{n+3}\,y_{n+3}+\bar\alpha_{n+4}\,y_{n+4}+
\cdots\bigg],
\end{align}
where
\[
\bar\alpha_{n+k}=\la^{n-2}(\la-1)^2\,\alpha_{n+k},
\]
and
\begin{equation}\label{baG}
|\bar\alpha_{n+k}|\leq |\la|^{n-1}\,\frac{(|\la|+1)^2}{(|\la|-1)^2}.
\end{equation}
By inserting (\ref{xG}) into (\ref{s2g}) and (\ref{u2G}), we have
\[
x=(\lambda I-G_n)^{-1}y=\frac{1}{q_{n+1}(\la)} \, (A+B)\,y
\]
where the matrix representations of $A$ and $B$ are
given by
\[
\textbf{A}
=
\left( \begin{array}{cccccccccc}
(\la-1)^2\la^{n-2} &     \rho  &     \rho\,\la  &\rho\,\la^2 &\cdots &\rho\,\la^{n-2}
&\la^{n-1}  &\bar\alpha_{n+2}  &\bar\alpha_{n+3}& \cdots      \\
(\la-1)^2\la^{n-3} &  \frac{\rho}{\la}  &    \rho  & \rho\,\la   &\cdots
&\rho\,\la^{n-3}  & \la^{n-2}  & \frac{\bar\alpha_{n+2}}{\la} &
\frac{\bar\alpha_{n+3}}{\la} & \cdots      \\
(\la-1)^2\la^{n-4} &  \frac{\rho}{\la^2}  &    \frac{\rho}{\la}  & \rho  &\cdots
&\rho\,\la^{n-4}  & \la^{n-3}  & \frac{\bar\alpha_{n+2}}{\la^2} &
\frac{\bar\alpha_{n+3}}{\la^2} & \cdots      \\
(\la-1)^2\la^{n-5} &  \frac{\rho}{\la^3}  &    \frac{\rho}{\la^2}  & \frac{\rho}{\la}  &\cdots
&\rho\,\la^{n-5}  & \la^{n-4}  & \frac{\bar\alpha_{n+2}}{\la^3} &
\frac{\bar\alpha_{n+3}}{\la^3} & \cdots      \\
 \vdots & \vdots   & \vdots      & \vdots   & \cdots      & \vdots
 & \vdots     & \vdots   & \vdots    &\cdots \\
(\la-1)^2\la &  \frac{\rho}{\la^{n-3}}  &  \frac{\rho}{\la^{n-4}}  &
\frac{\rho}{\la^{n-5}} & \cdots & \rho\,\la & \la^2 &
\frac{\bar\alpha_{n+2}}{\la^{n-3}} & 
\frac{\bar\alpha_{n+3}}{\la^{n-3}} & \cdots    \\
(\la-1)^2 & \frac{\rho}{\la^{n-2}}  &  \frac{\rho}{\la^{n-3}}  &
\frac{\rho}{\la^{n-4}} & \cdots & \rho & \la &
\frac{\bar\alpha_{n+2}}{\la^{n-2}} & 
\frac{\bar\alpha_{n+3}}{\la^{n-2}} & \cdots      \\
(\la-1)^2\,\frac{1}{\la} &  \frac{\rho}{\la^{n-1}}  &  \frac{\rho}{\la^{n-2}}  &
\frac{\rho}{\la^{n-3}} & \cdots & \frac{\rho}{\la} & 1 &
\frac{\bar\alpha_{n+2}}{\la^{n-1}} & 
\frac{\bar\alpha_{n+3}}{\la^{n-1}} &
\cdots     \\
(\la-1)^2\,\frac{\rho}{\la^2} &  \frac{\rho^2}{\la^{n}}  &  \frac{\rho^2}{\la^{n-1}}  &
\frac{\rho^2}{\la^{n-2}} & \cdots & \frac{\rho^2}{\la^2} & \frac{\rho}{\la} &
\frac{\rho\,\bar\alpha_{n+2}}{\la^{n}} & 
\frac{\rho\,\bar\alpha_{n+3}}{\la^{n}} &
\cdots     \\ 
(\la-1)^2\,\frac{2\,\rho}{\la^3} &  \frac{2\,\rho^2}{\la^{n+1}}  &
\frac{2\,\rho^2}{\la^{n}}  & 
\frac{2\,\rho^2}{\la^{n-1}} & \cdots & \frac{2\,\rho^2}{\la^3} & \frac{2\,\rho}{\la^2} &
\frac{2\,\rho\,\bar\alpha_{n+2}}{\la^{n+1}} & 
\frac{2\,\rho\,\bar\alpha_{n+3}}{\la^{n+1}} &
\cdots     \\ 
(\la-1)^2\,\frac{3\,\rho}{\la^4} &  \frac{3\,\rho^2}{\la^{n+2}}  &
\frac{3\,\rho^2}{\la^{n+1}}  & 
\frac{3\,\rho^2}{\la^{n}} & \cdots & \frac{3\,\rho^2}{\la^4} & \frac{3\,\rho}{\la^3} &
\frac{3\,\rho\,\bar\alpha_{n+2}}{\la^{n+2}} & 
\frac{3\,\rho\,\bar\alpha_{n+3}}{\la^{n+2}} &
\cdots     \\ 
(\la-1)^2\,\frac{4\,\rho}{\la^5} &  \frac{4\,\rho^2}{\la^{n+3}}  &
\frac{4\,\rho^2}{\la^{n+2}}  & 
\frac{4\,\rho^2}{\la^{n+1}} & \cdots & \frac{4\,\rho^2}{\la^5} & \frac{4\,\rho}{\la^4} &
\frac{4\,\rho\,\bar\alpha_{n+2}}{\la^{n+3}} & 
\frac{4\,\rho\,\bar\alpha_{n+3}}{\la^{n+3}} &
\cdots     \\ 
 \vdots & \vdots         & \vdots         & \vdots      & \cdots         & \vdots        & \vdots         & \vdots      & \vdots     &\cdots
\end{array}\right),
\]
and 
\[
\textbf{B}=
\left( \begin{array}{cccccccccccccc}
0 & 0 & 0 & 0 & 0 & \cdots &0 & 0 & 0 & 0  & 0 & 0 & 0& \cdots \\
0 & \frac{1}{\la} & 0 & 0 & 0 & \cdots & 0 & 0 & 0 & 0  & 0 & 0 & 0& \cdots \\
0 & \frac{1}{\la^2} & \frac{1}{\la} & 0 & 0 & \cdots & 0 & 0 & 0 & 0  & 0 & 0
& 0& \cdots\\ 
0 & \frac{1}{\la^3} & \frac{1}{\la^2} & \frac{1}{\la} & 0 & \cdots & 0 & 0 & 0
& 0  & 0 & 0 & 0& \cdots \\
\vdots & \vdots & \vdots & \vdots & \vdots & \cdots &\vdots & \vdots &\vdots &
\vdots & \vdots & \vdots & \vdots & \cdots \\
0 & \frac{1}{\la^{n-1}} & \frac{1}{\la^{n-2}} &
\frac{1}{\la^{n-3}} & \frac{1}{\la^{n-4}} & \cdots & \frac{1}{\la} & 0 & 0 & 0
 & 0 & 0 & 0 & \cdots \\
0 & \frac{\rho}{\la^{n}} & \frac{\rho}{\la^{n-1}} &
\frac{\rho}{\la^{n-2}} & \frac{\rho}{\la^{n-3}} & \cdots & \frac{\rho}{\la^2}
& \frac{1}{\la} & 0 & 0 
 & 0 & 0 & 0 & \cdots \\
0 & \frac{2\,\rho}{\la^{n+1}} & \frac{2\,\rho}{\la^{n}} &
\frac{2\,\rho}{\la^{n-1}} & \frac{2\,\rho}{\la^{n-2}} & \cdots & \frac{2\,\rho}{\la^3}
& \frac{2}{\la^2} & \frac{1}{\la} & 0 
 & 0 & 0 & 0 & \cdots \\
0 & \frac{3\,\rho}{\la^{n+2}} & \frac{3\,\rho}{\la^{n+1}} &
\frac{3\,\rho}{\la^{n}} & \frac{3\,\rho}{\la^{n-1}} & \cdots & \frac{3\,\rho}{\la^4}
& \frac{3}{\la^3} & \frac{3}{2\,\la^2} &  \frac{1}{\la}
 & 0 & 0 & 0 & \cdots \\
0 & \frac{4\,\rho}{\la^{n+3}} & \frac{4\,\rho}{\la^{n+2}} &
\frac{4\,\rho}{\la^{n+1}} & \frac{4\,\rho}{\la^{n}} & \cdots & \frac{4\,\rho}{\la^5}
& \frac{4}{\la^4} & \frac{4}{2\,\la^3} &  \frac{4}{3\,\la^2} & \frac{1}{\la}
 & 0 & 0 & \cdots \\
0 & \frac{5\,\rho}{\la^{n+4}} & \frac{5\,\rho}{\la^{n+3}} &
\frac{5\,\rho}{\la^{n+2}} & \frac{5\,\rho}{\la^{n+1}} & \cdots & \frac{5\,\rho}{\la^6}
& \frac{5}{\la^5} & \frac{5}{2\,\la^4} &  \frac{5}{3\,\la^3} & \frac{5}{4\,\la^2}
& \frac{1}{\la} & 0 & \cdots \\
\vdots & \vdots & \vdots & \vdots & \vdots & \cdots &\vdots & \vdots &\vdots &
\vdots & \vdots & \vdots & \vdots & \cdots
\end{array}\right),
\]
respectively. For $|\la|>1$, by using differentiation of the geometric series,
(\ref{baG}), and the argument used in (\ref{aG}), we have $\|A\|_1\leq \infty$ and 
$\|B\|_1\leq \infty$. Thus, for $|\la|>1$ and  
$q_{n+1}(\la)\neq 0$, the operator $\la I-G_n$ has a bounded inverse and its
resolvent set is given by 
(\ref{e3g}).

{\em Step 3.}
The point spectrum of the transposed operator $G_n^t$
consists of all $\lambda\in\mathbb{R}$, $|\la|\leq 1$,  such that
\[
(\lambda I - G_n^t)\,x=0, \qquad x\neq 0,\quad \|x\|_{\infty}< \infty.
\]
This is equivalent to
\begin{align*}
x_k&=\la^{k-1}\, x_1, \qquad k=2,\cdots,n, \\
x_{n+1}&=\frac{\la^n}{\rho} \, x_1, \\
x_{n+k}&=\frac{1}{k}\, \bigg( \frac{\la^{n+k-1}}{\rho}
  -\la^{k-2}-2\,\la^{k-3}-3\,\la^{k-4}-\cdots-(k-2)\,\la -(k-1)\bigg)\,x_1,\quad
k\geq 2.
\end{align*}
If $\la\neq 1$, we have
\begin{align*}
x_{n+k}&=\frac{1}{k}\, \bigg[ \frac{\la^{n+k-1}}{\rho}
  -\big(\la^{k-2}+2\,\la^{k-3}+\cdots+(k-2)\,\la
  +(k-1)\big)\frac{(\la-1)^2}{(\la-1)^2}\bigg]\,x_1 \\
  &=\frac{1}{k}\,\bigg(\frac{\la^{n+k-1}}{\rho} - \frac{\la^k-k\,\la
    +(k-1)}{(\la-1)^2}\bigg)\, x_1,\quad k\geq 2.
\end{align*}
Therefore
\[
|x_{n+k}|\leq \frac{1}{k}\bigg( \frac{1}{\rho}+\frac{1+k+k-1}{|\la
  -1|^2}\bigg)\,x_1 \leq \bigg(\frac{1}{2\,\rho}+\frac{2}{|\la-1|^2}\bigg)\,x_1,\qquad k\geq 2,
\]
which implies $\|x\|_{\infty}<\infty$. For $\la =1$ we have
\begin{align*}
x_{n+k}&=\frac{1}{k}\, \bigg(\frac{1}{\rho} -1-2-3-\cdots
-(k-2)-(k-1)\bigg)\,x_1\\ 
&=\frac{1}{k}\, \bigg(\frac{1}{\rho}
-\frac{(k-1)\,k}{2}\bigg)\,x_1 ,
\qquad k\geq 2,
\end{align*}
so $\|x\|_\infty = \infty$. We conclude that the point spectrum of
$G_n^t$ is given by (\ref{e4G}). This, in turn, implies (\ref{e5G}) and
(\ref{e5aG}) as described before.

\begin{example}\label{ex2}
The lesion forming plant pathogen potato late blight 
(phytophthora infestans) 
grows radially on a leaf with a constant daily rate.
The latency period for a lesion to become infectious is five days, and the
sporulating area is infectious for one day. 
In \cite{PSV05} the epidemic spread of such 
pathogen is modeled with the infinite dimensional Leslie matrix of the form of
$\Gamma_5$ as defined in (\ref{HdefGa}). 
Further, the upper bound for the speed
of invasion in computed via minimization of the largest eigenvalue
$\la_{\max}(\Gamma_5(s))$. 
From Remark \ref{rem2} it follows that this eigenvalue is the largest unique
positive root of $q_6(\lambda)$ from (\ref{polq}).
Here the parameter $\rho$ has the form 
$\rho(s)=\textrm{const}\times M(s)$ where $M(s)$ is some moment-generating function (for example,
$M(s)=\exp(\sigma^2s^2/2)$ for the Gaussian kernel or
$M(s)=1/(1-\sigma^2s^2)$ for the Laplace kernel).
Here $\Gamma_5(s)$ appears naturally due to the fact that the
considered pathogen has a latency period of five days.
It is interesting
that $\la_{\max}(\Gamma_5(s))$ can be computed analytically:
\[
\lambda_{\max}(\Gamma_5)= \frac{1}{3}+\frac{2^{1/3}}{3\left(2+27\sqrt{\rho}+
\sqrt{108\sqrt{\rho}+729\rho}\right)^{1/3}}+
\frac{\left(2+27\sqrt{\rho}+
\sqrt{108\sqrt{\rho}+729\rho}\right)^{1/3}}{3\cdot 2^{1/3}}.
\]
The speed of invasion is bounded by 
\[
v^* = \min_{0<s <\hat{s}} \frac1s \ln \left[\lambda_{\max}(\Gamma_5(s)) \right],
\label{Cdef}
\]
where $\hat{s}$ is the maximum $s$ for which $M(s)$ is
defined. For details about a rather complex
derivation of this model we refer the reader to \cite{PSV05}.
\end{example}



{\em Acknowledgment.} The author wishes to thank the anonymous referee for the
valuable comments and remarks, in particular for correcting an error in the
previous version of the manuscript and for providing Example \ref{ex1}.
The author also wishes to thank Jussi Behrndt from the Technische Universit\"at
Berlin for helpful discussions.

\end{document}